\documentclass[12pt]{article}
\usepackage{amsmath,amsthm, amssymb}
\usepackage{amssymb,latexsym}

\newtheorem{theorem}{Theorem}

\newtheorem{lemma}{Lemma}

\begin{document}

\baselineskip=17pt

\title{\bf A Bombieri -- Vinogradov type result for exponential sums over Piatetski-Shapiro primes}

\author{\bf S. I. Dimitrov}
\date{2022}
\maketitle
\begin{abstract}
In this paper, we establish a theorem of Bombieri -- Vinogradov type for exponential sums over Piatetski-Shapiro primes
$p= [n^{1/\gamma}]$ with $\frac{865}{886}<\gamma < 1$.
\\
\quad\\
\textbf{Keywords}:  Bombieri -- Vinogradov type result $\cdot$ Exponential sum $\cdot$ Piatetski-Shapiro Primes\\
\quad\\
{\bf  2020 Math.\ Subject Classification}: 11L07 $\cdot$  11L20  
\end{abstract}
\section{Introduction and statements of the result}
\indent

Let $\mathbb{P}$ denotes the set of all prime numbers.
In 1953 Piatetski-Shapiro \cite{Shapiro} has shown that for any fixed $\frac{11}{12}<\gamma< 1$ the set
\begin{equation*}
\mathbb{P}_\gamma=\{p\in\mathbb{P}\;\;|\;\; p= [n^{1/\gamma}]\;\; \mbox{ for some } n\in \mathbb{N}\}
\end{equation*}
is infinite.
The prime numbers of the form $p = [n^{1/\gamma}]$ are called Piatetski-Shapiro primes of type $\gamma$.
Denote
\begin{equation*}
P_\gamma(X)=\sum\limits_{p\leq X\atop{p=[n^{1/\gamma}]}}1\,.
\end{equation*}
Piatetski-Shapiro's result states that
\begin{equation*}
P_\gamma(X)\sim \frac{X^\gamma}{\log X}
\end{equation*}
for
\begin{equation*}
\frac{11}{12}<\gamma<1\,.
\end{equation*}
The best results up to now belongs to Rivat and Sargos \cite{Rivat-Sargos} with
\begin{equation}\label{Shapiroest}
P_\gamma(X)\sim \frac{X^\gamma}{\log X}
\end{equation}
for
\begin{equation*}
\frac{2426}{2817}<\gamma<1
\end{equation*}
and to Rivat and Wu \cite{Rivat-Wu} with
\begin{equation*}
P_\gamma(X)\gg\frac{X^\gamma}{\log X}
\end{equation*}
for
\begin{equation*}
\frac{205}{243}<\gamma<1.
\end{equation*}
On the other hand the celebrated Bombieri -- Vinogradov theorem is extremely important result in analytic number theory and has various applications.
It concerns the distribution of primes in arithmetic progressions, averaged over a range of moduli and states the following.
It asserts that when $A > 0$ is a fixed, then
\begin{equation*}
\sum\limits_{d\le \sqrt{X}/(\log X)^{A+5}}\max\limits_{y\le X}\max\limits_{(a,\,d)=1}
\Bigg|\sum_{p\le y\atop{p\equiv a\,( \textmd{mod}\, d)}} \, \log p-\frac{y}{\varphi(d)}\Bigg|\ll\frac{X}{\log^AX}\,,
\end{equation*}
where $\varphi (n)$ is Euler's function.

Recently  J. Li, M. Zhang and F. Xue \cite{Zhang-Li} obtained a mean value theorem of Bombieri-Vinigradov's type
for Piatetski-Shapiro primes. More precisely they showed that, when $\gamma$ is a real number satisfying
$\frac{85}{86}< \gamma < 1$ and $a\neq 0$ is a fixed integer, then for any given constant $A > 0$ and any sufficiently small $\varepsilon > 0$,
there holds
\begin{equation*}
\sum\limits_{d\le X^\theta\atop{(d, a)=1}}\Bigg|\sum\limits_{p\leq X\atop{p\equiv a\, (\textmd{mod}\, d)\atop{p=[n^{1/\gamma}]}}}1
-\frac{1}{\varphi(d)}P_\gamma(X)\Bigg|\ll\frac{X^\gamma}{\log^AX}\,,
\end{equation*}
where
\begin{equation*}
\theta=\frac{129}{4}\gamma-\frac{255}{8}-\varepsilon.
\end{equation*}
Weaker results were previously obtained from Peneva \cite{Peneva}, Wang and Cai\cite{Cai} and Lu \cite{Lu}.

Recently the author \cite{Dimitrov} established a new Bombieri -- Vinogradov type result for exponential sums over primes.
Let $1<c<3$, $c\neq2$, $0<\mu<1$ and $A>0$ be fixed. Then for $|t|<X^{\frac{1}{4}-c}$ the inequality
\begin{equation*}
\sum\limits_{d\le \sqrt{X}/(\log X)^{\frac{6A+34}{3}}}\max\limits_{y\le X}
\max\limits_{(a,\, d)=1}\Bigg|\sum\limits_{\mu y<p\leq y\atop{p\equiv a\, (\textmd{mod}\, d)}}e(t p^c)\log p
-\frac{1}{\varphi(d)}\int\limits_{\mu y}^{y}e(t x^c)\,dx\Bigg|\ll\frac{X}{\log^AX}
\end{equation*}
holds.

Motivated by these investigations in this paper we establish a new Bombieri -- Vinogradov type result
for exponential sums over Piatetski-Shapiro primes.
More precisely we establish the following theorem.
\begin{theorem}\label{Theorem}
Let $\frac{865}{886}<\gamma < 1$, $a\in\mathbb{Z}\setminus\{0\}$,
$1<c<3$, $c\neq2$ and $A>0$ be fixed.
Then for $|t|<X^{\frac{1}{4}-c}$ and any sufficiently small $\varepsilon > 0$ the inequality
\begin{equation*}
\sum\limits_{d\le X^\theta\atop{(d, a)=1}}\Bigg|\sum\limits_{p\leq X\atop{p\equiv a\, (\textmd{mod}\, d)\atop{p=[n^{1/\gamma}]}}}
p^{1-\gamma}e(t p^c) \log p-\frac{\gamma}{\varphi(d)}\int\limits_2^Xe(t y^c)\,dy\Bigg|\ll\frac{X}{\log^AX}
\end{equation*}
holds. Here
\begin{equation*}
\theta=\theta(\gamma)=\frac{443}{55}\gamma-\frac{173}{22}-\varepsilon\,.
\end{equation*}
\end{theorem}

\bigskip
\bigskip

\emph{Remark.}
From Theorem \ref{Theorem} it follows that 
\begin{equation*}
\theta(\gamma)\xrightarrow[\gamma\rightarrow1]{}\frac{21}{110}\,.
\end{equation*}

\section{Notations}
\indent

Assume that $X$ is a sufficiently large positive number and $k\geq8$ is a fixed natural number.
The letter $p$ with or without subscript will always denote prime numbers.
The notation $m\sim M$ means that $m$ runs through the interval $(M, 2M]$.
As usual $\varphi (n)$ is Euler's function,
$\tau(n)$ denotes the number of positive divisors of $n$ and $\Lambda(n)$ is von Mangoldt's function.
Let $(d,a)$ be the greatest common divisor of $d$ and $a$.
Instead of $m\equiv n\,\pmod {k}$ we write for simplicity $m\equiv n(k)$.
Further $[t]$, $\{t\}$ and $\|t\|$ denote the integer part of $t$, the fractional part of $t$
and the distance from $t$ to the nearest integer, respectively. Moreover $e(t)$=exp($2\pi it$) and $\psi(t)=\{t\}-\frac{1}{2}$.
The letter $\varepsilon$ denotes an arbitrary small positive number, not the same in all appearances.
Throughout this paper we suppose that  $\frac{865}{886}<\gamma < 1$. Denote
\begin{align}
\label{theta}
&\theta= \frac{443}{55}\gamma-\frac{173}{22}-\varepsilon\,;\\
\label{D}
&D=X^\theta\,.
\end{align}

\section{Preliminary lemmas}
\indent


\begin{lemma}\label{Bomb-Vin-Dim} Let $1<c<3$, $c\neq2$, $0<\mu<1$ and $A>0$ be fixed. Then for $|t|<X^{\frac{1}{4}-c}$ the inequality
\begin{equation*}
\sum\limits_{d\le \sqrt{X}/(\log X)^{\frac{6A+34}{3}}}\max\limits_{y\le X}
\max\limits_{(a,\, d)=1}\Bigg|\sum\limits_{\mu y<p\leq y\atop{p\equiv a\, ( d)}}e(t p^c)\log p
-\frac{1}{\varphi(d)}\int\limits_{\mu y}^{y}e(t x^c)\,dx\Bigg|\ll\frac{X}{\log^AX}
\end{equation*}
holds.
\end{lemma}
\begin{proof}
See (\cite{Dimitrov}, Lemma 18).
\end{proof}

\begin{lemma}\label{Exponentpairs1}
Let $|f^{(m)}(u)|\asymp YX^{1-m}$  for $1\leq X<u\leq X_1\leq2X$ and $m\geq1$.\\
Then
\begin{equation*}
\bigg|\sum_{X<n\le X_1}e(f(n))\bigg|\ll Y^\varkappa X^\lambda +Y^{-1},
\end{equation*}
where $(\varkappa, \lambda)$ is any exponent pair.
\end{lemma}
\begin{proof}
See (\cite{Graham-Kolesnik}, Ch. 3).
\end{proof}

\begin{lemma}\label{ExponentpairofBourgain}
For every $\varepsilon > 0$, the pair $\Big(\frac{13}{84}+\varepsilon, \frac{55}{84}+\varepsilon\Big)$ is an exponent pair.
\end{lemma}
\begin{proof}
See (\cite{Bourgain}, Theorem 6).
\end{proof}

\begin{lemma}\label{Heath-Brown} Let $G(n)$ is an arithmetic function.
Then
\begin{equation}\label{bilinearforms}
\bigg|\sum\limits_{n\sim N}\Lambda(n)G(n)\bigg|
\ll N^\varepsilon\max\bigg| \mathop{\sum\sum}_{ml\sim N\atop{m\sim M}}a(m)b(l)G(ml) \bigg|
\end{equation}
where the maximum is taken over all bilinear forms with coefficients satisfying one of the following three cases
\begin{align}
\label{Case1}
&|a(m)|\leq 1\,,\quad\quad  b(l)= 1\,,\\
\label{Case2}
&|a(m)|\leq 1\,,\quad\quad  b(l)= \log l\,,\\
\label{Case3}
&|a(m)|\leq 1\,,\quad \quad |b(l)|\leq 1
\end{align}
and also satisfying in all cases
\begin{equation}\label{MX}
M\leq N\,.
\end{equation}
\end{lemma}

\begin{proof}
See (\cite{Heath1}).
\end{proof}

The sums of cases \eqref{Case1}, \eqref{Case2} are called sums of Type I and are denoted by $S_I$.
The sums from case \eqref{Case3} are called sums of Type II and are denoted by $S_{II}$.

\begin{lemma}\label{Balog} Let $a$, $b$, $c$ be real numbers such that
\begin{align}
\label{0a1b}
&0 < a < 1\,,\quad b<\frac{2}{3}\,,\\
\label{0bc1}
&0 < b < c < 1 \,,\\
\label{1cb}
&1 - c < c - b\,,\\
\label{1ac}
&1 - a < \frac{c}{2} \,.
\end{align}
Then then \eqref{bilinearforms} still holds when \eqref{MX} is replaced by the conditions
\begin{align*}
&M\leq N^a\quad\quad  \mbox{for Type I sums}\,,\\
&N^b\leq M\leq N^c\quad\quad  \mbox{for Type II sums}\,.
\end{align*}
\end{lemma}

\begin{proof}
See (\cite{Balog}, Proposition 1).
\end{proof}

\begin{lemma}\label{thenumberofsolutions}
Let $\frac{1}{2}<\gamma<1$, $H \geq1$, $N \geq1$ and $\Delta>0$.
Denote by $\mathcal{N}(\Delta)$ the number of solutions of the inequality
\begin{equation*}
\big|h_1n_1^\gamma -h_2n_2^\gamma\big|\leq \Delta\,, \quad  h_1, h_2 \sim H\,,   \quad   n_1, n_2 \sim N\,.
\end{equation*}
Then
\begin{equation*}
\mathcal{N}(\Delta)\ll \Delta H N^{2-\gamma} +H N \log(H N)\,.
\end{equation*}
\end{lemma}
\begin{proof}
See (\cite{Heath2}).
\end{proof}

\section{Outline of the proof}
\indent

We recall  that
\begin{equation}\label{tlimits}
|t|\leq X^{\frac{1}{4}-c}\,.
\end{equation}
In order to prove Theorem \ref{Theorem} we need to establish the following estimates of exponential sums
\begin{align}
\label{Reduction1}
&\sum\limits_{d\le D\atop{(d, a)=1}}\Bigg|\sum\limits_{p\leq X\atop{p\equiv a\, (d)}}
p^{1-\gamma}\big((p+1)^\gamma-p^\gamma\big)e(t p^c) \log p-\frac{\gamma}{\varphi(d)}\int\limits_2^Xe(t y^c)\,dy\Bigg|\ll\frac{X}{\log^AX}\,,\\
\label{Reduction2}
&\sum\limits_{d\le D\atop{(d, a)=1}}\Bigg|\sum\limits_{p\leq X\atop{p\equiv a\, ( d)}}
p^{1-\gamma}\big(\psi(-(p+1)^\gamma)-\psi(-p^\gamma)\big)e(t p^c) \log p\Bigg|\ll\frac{X}{\log^AX}\,.
\end{align}
The estimate \eqref{Reduction1} follows immediately from Lemma \ref{Bomb-Vin-Dim}.
Thus it remains to prove  \eqref{Reduction2}.
Using the simplest splitting up argument, we obtain that to prove \eqref{Reduction2} it is sufficient to establish
\begin{equation}\label{Reduction3}
\sum\limits_{d\le D\atop{(d, a)=1}}\Bigg|\sum\limits_{n\sim N\atop{n\equiv a\, (d)}}
\Lambda(n)e(t n^c)\big(\psi(-(n+1)^\gamma)-\psi(-n^\gamma)\big) \Bigg|\ll\frac{X^\gamma}{\log^AX}
\end{equation}
for any $N\leq X$. Consider the case when
\begin{equation*}
N\leq X^{1-\varepsilon}\,.
\end{equation*}
We have
\begin{align*}
&\sum\limits_{d\le D\atop{(d, a)=1}}\Bigg|\sum\limits_{n\sim N\atop{n\equiv a\, (d)}}
\Lambda(n)e(t n^c)\big(\psi(-(n+1)^\gamma)-\psi(-n^\gamma)\big) \Bigg|\\
&\ll\sum\limits_{d\le D\atop{(d, a)=1}}\Bigg|\sum\limits_{n\sim N\atop{n\equiv a\, ( d)}}
\Lambda(n)e(t n^c)\big((n+1)^\gamma-n^\gamma\big) \Bigg|\\
&+\sum\limits_{d\le D\atop{(d, a)=1}}\Bigg|\sum\limits_{n\sim N\atop{n\equiv a\, ( d)}}
\Lambda(n)e(t n^c)\big([-n^\gamma]-[-(n+1)^\gamma]\big) \Bigg|\\
&\ll(\log X)\sum\limits_{n\sim N}n^{\gamma-1}\tau(n-a)+(\log X)\sum\limits_{n\sim N\atop{n=[m^{1/\gamma}]}}\tau(n-a)\\
&\ll N^{\gamma+\frac{\varepsilon}{2}} \ll\frac{X^\gamma}{\log^AX}\,.
\end{align*}
Henceforth we assume that
\begin{equation}\label{XNX}
X^{1-\varepsilon}\leq N\leq X\,.
\end{equation}
Obviously when
\begin{equation}\label{thetaepsilon}
\theta \leq \frac{1-\varepsilon}{2}
\end{equation}
then \eqref{D}, \eqref{XNX} and \eqref{thetaepsilon} give us
\begin{equation}\label{NDtheta}
N^\theta \leq D\leq N^{\theta+\frac{\varepsilon}{2}}\,.
\end{equation}
We recall the well-known expansions
\begin{equation}\label{expansion}
\psi(t)=-\sum\limits_{1\leq|h|\leq H}\frac{e(ht)}{2\pi i h}
+\mathcal{O}\Bigg(\min\left(1,\frac{1}{H\|t\|}\right)\Bigg)\,,
\end{equation}
\begin{equation*}
\min\left(1,\frac{1}{H\|t\|}\right)=\sum\limits_{k=-\infty}^{\infty}b_ke(k t)\,,
\end{equation*}
where
\begin{equation*}
|b_k|\ll\min\left(\frac{\log 2H}{k},\frac{1}{|k|}, \frac{H}{|k|^2}\right)\,.
\end{equation*}
Using \eqref{expansion} for the left-hand side of \eqref{Reduction3} we deduce
\begin{equation}\label{Reduction4}
\sum\limits_{d\le D\atop{(d, a)=1}}\Bigg|\sum\limits_{n\sim N\atop{n\equiv a\, (d)}}
\Lambda(n)e(t n^c)\big(\psi(-(n+1)^\gamma)-\psi(-n^\gamma)\big) \Bigg|\ll \Gamma_1+\Gamma_2+\Gamma_3\,,
\end{equation}
where
\begin{align}
\label{Gamma1}
&\Gamma_1=\sum\limits_{d\le D\atop{(d, a)=1}}\sum\limits_{1\leq h\leq H}\frac{1}{h}
\Bigg|\sum\limits_{n\sim N\atop{n\equiv a\, ( d)}}
\Lambda(n)e(t n^c)\Big(e(-hn^\gamma)-e(-h(n+1)^\gamma)\Big) \Bigg| \,, \\
\label{Gamma2}
&\Gamma_2=\sum\limits_{d\le D\atop{(d, a)=1}}\sum\limits_{n\sim N\atop{n\equiv a\, ( d)}}
\Lambda(n)\sum\limits_{k=-\infty}^{\infty}b_ke(k n^\gamma)\,,\\
\label{Gamma3}
&\Gamma_3=\sum\limits_{d\le D\atop{(d, a)=1}}\sum\limits_{n\sim N\atop{n\equiv a\, ( d)}}
\Lambda(n)\sum\limits_{k=-\infty}^{\infty}b_ke(k (n+1)^\gamma) \,.
\end{align}

We shall estimate $\Gamma_2$, $\Gamma_3$ and $\Gamma_1$, respectively,
in the sections \ref{SectionGamma3} and \ref{SectionGamma1}.
In section \ref{Sectionfinal} we shall finalize the proof of Theorem \ref{Theorem}.

\section{Upper bound of  $\mathbf{\Gamma_2}$ and $\mathbf{\Gamma_3}$}\label{SectionGamma3}
\indent

Arguing as in  \cite{Zhang-Li} for the sums denoted by \eqref{Gamma2} and \eqref{Gamma3} we have that
\begin{equation}\label{Gamma23est}
\Gamma_2, \Gamma_3\ll\frac{X^\gamma}{\log^AX}\,,
\end{equation}
provided that
\begin{equation}\label{Hgamma}
H=N^{1-\gamma+\varepsilon}\,, \quad \gamma> \frac{1}{2}+\theta\,.
\end{equation}

\section{Upper bound of  $\mathbf{\Gamma_1}$}\label{SectionGamma1}
\indent

Now we shall estimate $\Gamma_1$. Proceeding as in \cite{Zhang-Li} for the sum denoted by \eqref{Gamma1} we obtain
\begin{equation}\label{Gamma1est1}
\Gamma_1\ll X^{\gamma-1}\sum\limits_{n\sim N}\Lambda(n) G(n)\,,
\end{equation}
where
\begin{align}
\label{Gn}
&G(n)=\sum\limits_{1\leq h\leq H}F_h(n)e(tn^c-hn^\gamma)\,,\\
\label{Fhn}
&F_h(n)=\sum\limits_{d\le D\atop{(d, a)=1\atop{d|n- a}}}c(d,h)\,,\\
\label{cdhn}
&|c(d,h)|=1\,.
\end{align}
Therefore, it remains to show that
\begin{equation}\label{LambdaGN}
\bigg|\sum\limits_{n\sim N}\Lambda(n) G(n)\bigg|\ll \frac{N}{\log^AN}\,.
\end{equation}
In order to apply Lemma \ref{Balog} we need to find the upper bounds of the sums of Type I
and of Type II. As we mentioned, $S_I$ denotes sum of Type I and $S_{II}$ denotes sum of Type II.

\begin{lemma}\label{Exponentpairs2}
Let $1 \leq d \leq X$ , $X < X_1\leq 2X$. Then
\begin{align*}
&\sum_{X<n\le X_1\atop{n\equiv a\, ( d)}}e(h_1 n^c+h_2 n^\gamma)\\
&\ll \min \left(\frac{X}{d}\,,  \, \big(|h_1|dX^{c-1} +  |h_2|dX^{\gamma-1} \big )^{-1}
+d^{\varkappa-\lambda}|h_1|^\varkappa X^{\varkappa c-\varkappa+\lambda}
+d^{\varkappa-\lambda}|h_2|^\varkappa X^{\varkappa\gamma-\varkappa+\lambda}\right)\,,
\end{align*}
where $(\varkappa, \lambda)$ is any exponent pair.
\end{lemma}
\begin{proof}
Let $b$ is an integer such that $1 \leq b \leq d$ and $b\equiv a \, ( d)$. Then
\begin{equation}\label{Exponentpairs2est1}
\sum_{X<n\le X_1\atop{n\equiv a\, ( d)}}e(h_1 n^c+h_2 n^\gamma)
=\sum_{\frac{X-b}{d}<m\le \frac{X_1-b}{d}}e\Big(h_1 (b+md)^c+h_2 (b+md)^\gamma\Big)\,.
\end{equation}
Denote
\begin{equation*}
f(y)=h_1 (b+yd)^c+h_2 (b+yd)^\gamma\,.
\end{equation*}
We have
\begin{equation}\label{derivative}
|f^{(m)}(y)|\asymp \big(|h_1|dX^{c-1} +  |h_2|dX^{\gamma-1} \big ) \left(\frac{X}{d}\right)^{1-m}
\end{equation}
for
\begin{equation*}
X<y\leq X_1 \quad \mbox{ and } \quad  m\geq1\,.
\end{equation*}
Using \eqref{Exponentpairs2est1}, \eqref{derivative}, the trivial estimation and Lemma \ref{Exponentpairs1}
by any exponent pair $(\varkappa, \lambda)$ we establish the statement in the lemma.
\end{proof}

\begin{lemma}\label{SIest} Assume that
\begin{equation}\label{Conditions1}
M\ll N^{\frac{31}{3}\gamma-\frac{19}{2}-\theta-\varepsilon} \,.
\end{equation}
Then
\begin{equation*}
S_I\ll N^{1-\varepsilon}\,.
\end{equation*}
\end{lemma}
\begin{proof}

Splitting the range of $l$  into dyadic subintervals of the form $(L, 2L]$ so that $ML\sim N$,
changing the order of summation, using \eqref{tlimits}, \eqref{XNX}, \eqref{NDtheta}, \eqref{Hgamma},
\eqref{Conditions1} and Lemma \ref{Exponentpairs2}
with the exponent pair
\begin{equation*}
A^4B(0, 1)=\left(\frac{1}{62},\frac{57}{62}\right)
\end{equation*}
we deduce
\begin{align}\label{SIest1}
S_{I}&\ll(\log X)^2\sum\limits_{1\leq h\leq H}\sum\limits_{d\le D\atop{(d, a)=1}}c(d,h)
\sum\limits_{m\sim M}a(m) \sum\limits_{l\sim L\atop{ml\sim N \atop{ml\equiv a\,(d)}}}e\big(t(ml)^c-h(ml)^\gamma\big) \nonumber\\
&\ll (\log X)^2\sum\limits_{1\leq h\leq H}\sum\limits_{d\le D\atop{(d, a)=1}}
\sum\limits_{m\sim M}\bigg|\sum\limits_{l\sim L\atop{ml\sim N \atop{ml\equiv a\,(d)}}}e\big(t(ml)^c-h(ml)^\gamma\big)\bigg|\nonumber\\
&\ll (\log X)^2\sum\limits_{1\leq h\leq H}\sum\limits_{d\le D\atop{(d, a)=1}}
\sum\limits_{m\sim M}\Big( \big(|t|dMN^{c-1} +  hdMN^{\gamma-1} \big )^{-1}\nonumber\\
&+d^{-\frac{28}{31}}|t|^\frac{1}{62}M^{-\frac{28}{31}} N^\frac{c+56}{62}
+d^{-\frac{28}{31}}h^\frac{1}{62}M^{-\frac{28}{31}}N^\frac{\gamma+56}{62}\Big)\nonumber\\
&\ll (\log X)^2\sum\limits_{1\leq h\leq H}\sum\limits_{d\le D\atop{(d, a)=1}}
\sum\limits_{m\sim M}\Big( h^{-1}d^{-1}M^{-1}N^{1-\gamma} \nonumber\\
&+d^{-\frac{28}{31}}N^{\left(\frac{1}{4}-c\right)\frac{1}{62}}M^{-\frac{28}{31}} N^\frac{c+56}{62}
+d^{-\frac{28}{31}}h^\frac{1}{62}M^{-\frac{28}{31}}N^\frac{\gamma+56}{62}\Big)\nonumber\\
&\ll (\log X)^4\Big( N^{1-\gamma} +H M^\frac{3}{31}N^\frac{225}{248}D^\frac{3}{31}
+H^\frac{63}{62} M^\frac{3}{31}N^\frac{\gamma+56}{62}D^\frac{3}{31}\Big)\nonumber\\
&\ll (\log X)^4\Big( N^{1-\gamma}+H^\frac{63}{62} M^\frac{3}{31}N^\frac{\gamma+56}{62}D^\frac{3}{31}\Big)\nonumber\\
&\ll N^{1-\varepsilon}\,.
\end{align}
Bearing in mind \eqref{Conditions1} and \eqref{SIest1} we establish the statement in the lemma.
\end{proof}

\begin{lemma}\label{SIIest} Assume that
\begin{align}
\label{Conditions2a}
&0<\theta \leq \frac{21}{110}-\varepsilon\,,\\
\label{Conditions2b}
&\frac{865}{886}+ \frac{55}{443}\theta+\varepsilon <\gamma <1\,,\\
\label{Conditions2c}
&N^{\frac{313}{44}-\frac{388}{55}\gamma+\theta+\varepsilon} \ll M\ll N^{\gamma-\varepsilon} \,.
\end{align}
Then
\begin{equation*}
S_{II}\ll N^{1-\varepsilon}\,.
\end{equation*}
\end{lemma}
\begin{proof}

Splitting the range of $l$ and $h$ into dyadic subintervals of the form $(L, 2L]$ and $(K, 2K]$
so that $ML\sim N$ and
\begin{equation}\label{KH}
\frac{1}{2}\leq K\leq H
\end{equation}
we obtain
\begin{equation*}
S_{II}\ll(\log X)^2\sum\limits_{m\sim M}\bigg| \sum\limits_{l\sim L\atop{ml\sim N}}
\sum\limits_{h\sim K}b(l)F_h(ml)e\big(t(ml)^c-h(ml)^\gamma\big) \bigg|\,.
\end{equation*}
It is easy to see that
\begin{equation*}
0<hl^\gamma\leq4KL^\gamma\,.
\end{equation*}
Let $T$ be a parameter, which will be determined later.
We decompose the pairs $(h, l)$ into sets  $\mathfrak{S}_y$  $(1\leq y \leq T)$  defined by
\begin{equation*}
\mathfrak{S}_y=\left\{(h, l)\;\; |\quad  h\sim K\,, \quad l\sim L\,, \quad  \frac{4KL^\gamma (y-1)}{T}<hl^\gamma\leq\frac{4KL^\gamma y}{T} \right \}\,.
\end{equation*}
Hence
\begin{equation*}
S_{II}\ll(\log X)^2\sum\limits_{1\leq y \leq T}\sum\limits_{m\sim M}\bigg| \mathop{\sum\sum}_{(h,l)\in\mathfrak{S}_y\atop{ml\sim N}}
b(l)F_h(ml)e\big(t(ml)^c-h(ml)^\gamma\big) \bigg|\,.
\end{equation*}
Using Cauchy's inequality we get
\begin{align}\label{SIIest1}
|S_{II}|^2&\ll(\log X)^4TM\sum\limits_{1\leq y \leq T}\sum\limits_{m\sim M}\bigg| \mathop{\sum\sum}_{(h,l)\in\mathfrak{S}_y\atop{ml\sim N}}
b(l)F_h(ml)e\big(t(ml)^c-h(ml)^\gamma\big)\bigg|^2\nonumber\\
&\ll(\log X)^4TM\sum\limits_{1\leq y \leq T}\mathop{\sum\sum}_{(h_1,l_1)\in\mathfrak{S}_y}\mathop{\sum\sum}_{(h_2,l_2)\in\mathfrak{S}_y}
\bigg| \sum\limits_{m\sim M\atop{ml_1\sim N\atop{ml_2\sim N}}}F_{h_1}(ml_1)F_{h_2}(ml_2)\nonumber\\
&\hspace{60mm}\times e\Big(\big(l_1^c -l_2^c\big)tm^c-\big(h_1l_1^\gamma -h_2l_2^\gamma\big) m^\gamma\Big) \bigg|\nonumber\\
&\ll(\log X)^4TM
\mathop{\sum\limits_{h_1\sim K}\sum\limits_{h_2\sim K}\sum\limits_{l_1\sim L}\sum\limits_{l_2\sim L}}_{|\alpha|\leq4KL^\gamma T^{-1}}
\big| \Theta(t) \big|\,,
\end{align}
where
\begin{equation}\label{Theta}
\Theta(t)= \sum\limits_{m\sim M\atop{ml_1\sim N\atop{ml_2\sim N}}}F_{h_1}(ml_1)F_{h_2}(ml_2)
e\Big(\big(l_1^c -l_2^c\big)tm^c-\alpha m^\gamma\Big)\,,
\end{equation}
\begin{equation}\label{alpha}
\alpha=h_1l_1^\gamma -h_2l_2^\gamma\,.
\end{equation}
If the system of congruences
\begin{equation}\label{System}
\left| \begin{array}{c}
l_1m\equiv a\,(d_1)\\
l_2m\equiv a\,(d_2)
\end{array}
\right.
\end{equation}
has no solution then $\Theta(t)=0$. Assume that the system \eqref{System} has a solution.
Then there exists an integer $f=f(l_1, l_2, a, d_1, d_2)$ such that $(f,[d_1, d_2])=1$
and \eqref{System} is equivalent to $m\equiv f\,([d_1, d_2])$.
From  \eqref{Fhn}, \eqref{cdhn}, \eqref{Theta} and the last consideration it follows
\begin{align}\label{Thetaest1}
\Theta(t)&= \sum\limits_{d_1\le D\atop{(d_1, a)=1}}c(d_1,h_1) \sum\limits_{d_2\le D\atop{(d_2, a)=1}}c(d_2, h_2)
\sum\limits_{m\sim M\atop{ml_1\sim N, \, ml_2\sim N\atop{m\equiv f\,([d_1, d_2])}}}
e\Big(\big(l_1^c -l_2^c\big)tm^c-\alpha m^\gamma\Big)\nonumber\\
&\ll\sum\limits_{d_1\le D\atop{(d_1, a)=1}} \sum\limits_{d_2\le D\atop{(d_2, a)=1}}
\bigg|\sum\limits_{m\sim M\atop{ml_1\sim N, \, ml_2\sim N\atop{m\equiv f\,([d_1, d_2])}}}
e\Big(\big(l_1^c -l_2^c\big)tm^c-\alpha m^\gamma\Big)\bigg|\,.
\end{align}
According to Lemma \ref{ExponentpairofBourgain}, the pair
\begin{equation*}
A\left(\frac{13}{84}+\varepsilon, \frac{55}{84}+\varepsilon\right)=\left(\frac{13}{194}+\varepsilon,\frac{152}{194}+\varepsilon\right)
\end{equation*}
is an exponent pair.
Now \eqref{Thetaest1} and Lemma \ref{Exponentpairs2} with exponent pair
\begin{equation*}
\left(\frac{13}{194}+\varepsilon,\frac{152}{194}+\varepsilon\right)
\end{equation*}
imply
\begin{equation}\label{Thetaest2}
\Theta(t)\ll\sum\limits_{d_1\le D\atop{(d_1, a)=1}} \sum\limits_{d_2\le D\atop{(d_2, a)=1}}
\min \left(\frac{M}{[d_1, d_2]},  \, \Omega\right)\,,
\end{equation}
where
\begin{align}\label{Omega}
\Omega=&\big(M^{c-1}[d_1, d_2] |(l_1^c -l_2^c\big)t| + M^{\gamma-1} [d_1, d_2] |\alpha| \big )^{-1}\nonumber\\
&+[d_1, d_2]^{-\frac{139}{194}}|(l_1^c -l_2^c\big)t|^\frac{13}{194} M^{\frac{13c}{194}+\frac{139}{194}}
+[d_1, d_2]^{-\frac{139}{194}}|\alpha|^\frac{13}{194} M^{\frac{13\gamma}{194}+\frac{139}{194}}\,.
\end{align}
Denote
\begin{align}
\label{E1}
&E_1=\mathop{\sum\limits_{h_1\sim K}\sum\limits_{h_2\sim K}\sum\limits_{l_1\sim L}\sum\limits_{l_2\sim L}}_{|\alpha|\leq4KL^\gamma T^{-1}}
\sum\limits_{d_1\le D\atop{(d_1, a)=1}} \sum\limits_{d_2\le D\atop{(d_2, a)=1}}\frac{M}{[d_1, d_2]}\,,\\
\label{E2}
&E_2=\mathop{\sum\limits_{h_1\sim K}\sum\limits_{h_2\sim K}\sum\limits_{l_1\sim L}\sum\limits_{l_2\sim L}}_{|\alpha|\leq4KL^\gamma T^{-1}}
\sum\limits_{d_1\le D\atop{(d_1, a)=1}} \sum\limits_{d_2\le D\atop{(d_2, a)=1}}
\big(M^{c-1}[d_1, d_2] |(l_1^c -l_2^c\big)t| + M^{\gamma-1} [d_1, d_2] |\alpha| \big )^{-1}\,,\\
\label{E3}
&E_3=\mathop{\sum\limits_{h_1\sim K}\sum\limits_{h_2\sim K}\sum\limits_{l_1\sim L}\sum\limits_{l_2\sim L}}_{|\alpha|\leq4KL^\gamma T^{-1}}
\sum\limits_{d_1\le D\atop{(d_1, a)=1}} \sum\limits_{d_2\le D\atop{(d_2, a)=1}}
[d_1, d_2]^{-\frac{139}{194}}|(l_1^c -l_2^c\big)t|^\frac{13}{194} M^{\frac{13c}{194}+\frac{139}{194}}\,,\\
\label{E4}
&E_4=\mathop{\sum\limits_{h_1\sim K}\sum\limits_{h_2\sim K}\sum\limits_{l_1\sim L}\sum\limits_{l_2\sim L}}_{|\alpha|\leq4KL^\gamma T^{-1}}
\sum\limits_{d_1\le D\atop{(d_1, a)=1}} \sum\limits_{d_2\le D\atop{(d_2, a)=1}}
[d_1, d_2]^{-\frac{139}{194}}|\alpha|^\frac{13}{194} M^{\frac{13\gamma}{194}+\frac{139}{194}}\,.
\end{align}
If
\begin{equation*}
|\alpha|\leq M^{-\gamma}
\end{equation*}
then
\begin{equation*}
\frac{M}{[d_1, d_2]}\leq\frac{M^{1-\gamma}}{[d_1, d_2] |\alpha|}
\end{equation*}
and the contribution of $M [d_1, d_2]^{-1}$ to $|S_{II}|^2$ is
\begin{equation}\label{E'1est}
E'_1\ll(\log X)^7 TM^2 \mathcal{N}(M^{-\gamma})\,,
\end{equation}
where the elementary estimate
\begin{equation}\label{elementaryestimate}
\sum\limits_{d_1\le D} \sum\limits_{d_2\le D}\frac{1}{[d_1, d_2]}\ll (\log D)^3\,.
\end{equation}
is used.\\
If
\begin{equation*}
|\alpha|\geq M^{-\gamma}
\end{equation*}
then
\begin{equation*}
\frac{M}{[d_1, d_2]}\geq\frac{M^{1-\gamma}}{[d_1, d_2] |\alpha|}
\end{equation*}
and the contribution of $\big(M^{c-1}[d_1, d_2] |(l_1^c -l_2^c\big)t| + M^{\gamma-1} [d_1, d_2] |\alpha| \big )^{-1}$ to $|S_{II}|^2$ is
\begin{equation}\label{E'2est}
E'_2\ll(\log X)^8 TM^{2-\gamma} \max\limits_{M^{-\gamma}\leq\Delta\leq4KL^\gamma T^{-1}}\frac{\mathcal{N}(\Delta)}{\Delta}\,,
\end{equation}
where the estimate \eqref{elementaryestimate} has been used again.\\
Further the total contribution of the term $[d_1, d_2]^{-\frac{139}{194}}|(l_1^c -l_2^c\big)t|^\frac{13}{194} M^{\frac{13c}{194}+\frac{139}{194}}$
to $|S_{II}|^2$ is
\begin{align}\label{E'3est}
E'_3&\ll(\log X)^4TM L^{\frac{13c}{194}} M^{\frac{13c}{194}+\frac{139}{194}}|t|^{\frac{13}{194}}
\sum\limits_{d_1\le D} \sum\limits_{d_2\le D}\frac{1}{[d_1, d_2]^\frac{139}{194}}\mathcal{N}(4KL^\gamma T^{-1})\nonumber\\
&\ll (\log X)^4T L^{\frac{13c}{194}} M^{\frac{13c}{194}+\frac{333}{194}}N^{\left(\frac{1}{4}-c\right)\frac{13}{194}}
D^\frac{55}{194}\mathcal{N}(4KL^\gamma T^{-1})
\end{align}
where \eqref{tlimits}, \eqref{XNX} and the upper bound
\begin{align}\label{upperbound}
\sum\limits_{d_1\le D} \sum\limits_{d_2\le D}\frac{1}{[d_1, d_2]^\frac{139}{194}}
&\ll\sum\limits_{d_1\le D} \sum\limits_{d_2\le D}\left(\frac{(d_1, d_2)}{d_1 d_2}\right)^\frac{139}{194}
=\sum\limits_{1\leq r\le D} \sum\limits_{k_1\le \frac{D}{r}}
\sum\limits_{k_2\le \frac{D}{r}}\frac{1}{r^\frac{139}{194} k_1^\frac{139}{194} k_2^\frac{139}{194}}\nonumber\\
&\ll\sum\limits_{1\leq r\le D}\frac{1}{r^\frac{139}{194}}\left(\sum\limits_{k\le \frac{D}{r}}\frac{1}{k^\frac{139}{194}}\right)^2
\ll \sum\limits_{1\leq r\le D}\frac{1}{r^\frac{139}{194}}\left( \frac{D}{r}\right)^\frac{110}{194}\ll D^\frac{55}{194}\,.
\end{align}
are used.\\
It remains to note that the total contribution of the term
$[d_1, d_2]^{-\frac{139}{194}}|\alpha|^\frac{13}{194} M^{\frac{13\gamma}{194}+\frac{139}{194}}$ to $|S_{II}|^2$ is
\begin{align}\label{E'4est}
E'_4&\ll(\log X)^4TM M^{\frac{13\gamma}{194}+\frac{139}{194}}|\alpha|^{\frac{13}{194}}
\sum\limits_{d_1\le D} \sum\limits_{d_2\le D}\frac{1}{[d_1, d_2]^\frac{139}{194}}\mathcal{N}(4KL^\gamma T^{-1})\nonumber\\
&\ll (\log X)^4T M^{\frac{13\gamma}{194}+\frac{333}{194}}(KL^\gamma T^{-1})^{\frac{13}{194}} D^{\frac{55}{194}}\mathcal{N}(4KL^\gamma T^{-1})\,,
\end{align}
where the upper bound \eqref{upperbound} has been used again.

Bearing in mind  \eqref{SIIest1}, \eqref{Thetaest2} -- \eqref{E'1est}, \eqref{E'2est},  \eqref{E'3est} and \eqref{E'4est}
and Lemma \ref{thenumberofsolutions} we deduce
\begin{align}\label{SIIest2}
|S_{II}|^2&\ll(\log X)^7 TM^2 \mathcal{N}(M^{-\gamma})+(\log X)^8 TM^{2-\gamma} \max\limits_{M^{-\gamma}\leq\Delta\leq4KL^\gamma T^{-1}}\frac{\mathcal{N}(\Delta)}{\Delta}\nonumber\\
&+(\log X)^4T L^{\frac{13c}{194}} M^{\frac{13c}{194}+\frac{333}{194}}N^{\left(\frac{1}{4}-c\right)\frac{13}{194}}
D^\frac{55}{194}\mathcal{N}(4KL^\gamma T^{-1})\nonumber\\
&+(\log X)^4T M^{\frac{13\gamma}{194}+\frac{333}{194}}(KL^\gamma T^{-1})^{\frac{13}{194}} D^{\frac{55}{194}}\mathcal{N}(4KL^\gamma T^{-1})\nonumber\\
&\ll(\log X)^8 TM^{2-\gamma} \max\limits_{M^{-\gamma}\leq\Delta\leq4KL^\gamma T^{-1}}\Big( K L^{2-\gamma} + \Delta^{-1} K L \log X\Big)\nonumber\\
&+(\log X)^4\Big(T N^{\frac{13}{776}} M^{\frac{333}{194}} D^{\frac{55}{194}}
+T^{\frac{181}{194}}M^\frac{333}{194}D^{\frac{55}{194}}K^{\frac{13}{194}}N^{\frac{13\gamma}{194}}\Big)\Big(K^2 L^2T^{-1}+K L\log X\Big)\nonumber\\
&\ll(\log X)^9 \Big( TK N^{2-\gamma} + T K NM +K^2N^{\frac{1565}{776}} M^{-\frac{55}{194}} D^{\frac{55}{194}}
+ T KN^{\frac{789}{776}} M^{\frac{139}{194}} D^{\frac{55}{194}}\nonumber\\
&+ T^{-\frac{13}{194}} K^{\frac{401}{194}} N^{\frac{13\gamma}{194}+2}M^{-\frac{55}{194}} D^{\frac{55}{194}}
+ T^{\frac{181}{194}} K^{\frac{207}{194}} N^{\frac{13\gamma}{194}+1}M^{\frac{139}{194}} D^{\frac{55}{194}}\Big)\,.
\end{align}
We take $T$ such that
\begin{equation*}
 T K NM= T^{-\frac{13}{194}} K^{\frac{401}{194}} N^{\frac{13\gamma}{194}+2}M^{-\frac{55}{194}} D^{\frac{55}{194}}\,.
\end{equation*}
Therefore
\begin{equation}\label{T}
T=\big[  N^{\frac{13\gamma+194}{207}}M^{-\frac{83}{69}}K D^{\frac{55}{207}}\big]+1\,.
\end{equation}
Now  \eqref{NDtheta}, \eqref{Hgamma}, \eqref{Conditions2a} -- \eqref{KH}, \eqref{SIIest2} and \eqref{T} yield
\begin{align}\label{SIIest3}
|S_{II}|&\ll N^\varepsilon\Big(K^2 N^{\frac{608-194\gamma}{207}} M^{-\frac{83}{69}} D^{\frac{55}{207}}
+ K^2 N^{\frac{13\gamma+401}{207}} M^{-\frac{14}{69}} D^{\frac{55}{207}}+K^2N^{\frac{1565}{776}} M^{-\frac{55}{194}} D^{\frac{55}{194}}\nonumber\\
&+ K^2N^{\frac{10088\gamma+313867}{160632}} M^{-\frac{6511}{13386}} D^{\frac{22055}{40158}}
+ K^2 N^{\frac{26\gamma+388}{207}}M^{-\frac{28}{69}} D^{\frac{110}{207}}+K N^{2-\gamma}+K NM\nonumber\\
&+ KN^{\frac{789}{776}} M^{\frac{139}{194}} D^{\frac{55}{194}}
+K^{\frac{207}{194}} N^{\frac{13\gamma}{194}+1}M^{\frac{139}{194}} D^{\frac{55}{194}}\Big)^{\frac{1}{2}}\nonumber\\
&\ll N^\varepsilon\Big( N^{\frac{1022-608\gamma}{207}} M^{-\frac{83}{69}} D^{\frac{55}{207}}
+ N^\frac{815-401\gamma}{207} M^{-\frac{14}{69}} D^{\frac{55}{207}}+N^\frac{3117-1552\gamma}{776} M^{-\frac{55}{194}} D^{\frac{55}{194}}\nonumber\\
&+ N^\frac{635131-311176\gamma}{160632} M^{-\frac{6511}{13386}} D^{\frac{22055}{40158}}
+  N^\frac{802-388\gamma}{207}M^{-\frac{28}{69}} D^{\frac{110}{207}}+N^{3-2\gamma}+ N^{2-\gamma}M\nonumber\\
&+N^{\frac{1565-776\gamma}{776}} M^{\frac{139}{194}} D^{\frac{55}{194}}
+ N^\frac{401-194\gamma}{194}M^{\frac{139}{194}} D^{\frac{55}{194}}\Big)^{\frac{1}{2}}\nonumber\\
&\ll N^{1-\varepsilon}\,.
\end{align}
The lemma is proved.
\end{proof}
We are now in a good position to establish the validity of \eqref{LambdaGN}.

In Lemma \ref{Balog} we take
\begin{align*}
&a= \frac{31}{3}\gamma-\frac{19}{2}-\theta-\varepsilon\,, \\
&b=\frac{313}{44}-\frac{388}{55}\gamma+\theta+\varepsilon\,, \\
&c=\gamma-\varepsilon\,.
\end{align*}
The direct verification assures us that the inequalities \eqref{0a1b} --  \eqref{1ac},
\eqref{Hgamma} and \eqref{Conditions2b} are fulfilled.
From \eqref{theta}, Lemma \ref{Balog}, Lemma \ref{SIest} and Lemma \ref{SIIest} it follows that \eqref{LambdaGN} holds.

Bearing in mind \eqref{XNX}, \eqref{Gamma1est1} and \eqref{LambdaGN} we deduce
\begin{equation}\label{Gamma1est}
\Gamma_1\ll\frac{X^\gamma}{\log^AX}\,.
\end{equation}

\section{The end of the proof}\label{Sectionfinal}
\indent

Summarizing  \eqref{Reduction4}, \eqref{Gamma23est} and \eqref{Gamma1est} we establish Theorem \ref{Theorem}.

\vskip20pt
\footnotesize
\begin{flushleft}
S. I. Dimitrov\\
Faculty of Applied Mathematics and Informatics\\
Technical University of Sofia \\
8, St.Kliment Ohridski Blvd. \\
1756 Sofia, BULGARIA\\
e-mail: sdimitrov@tu-sofia.bg\\
\end{flushleft}

\end{document}